\documentstyle[11pt,twoside]{article}
\title{Asymptotic evaluation of an integral arising in quantum harmonic oscillator tunnelling probabilities}
\author{\sc R. B.\ Paris \\
{\em School of Engineering, Computing and Applied Mathematics}, \\
{\em University of Abertay Dundee, Dundee DD1 1HG, UK}
}
\begin{document}
\def\f#1#2{\mbox{${\textstyle \frac{#1}{#2}}$}}
\def\dfrac#1#2{\displaystyle{\frac{#1}{#2}}}
\def\boldal{\mbox{\boldmath $\alpha$}}
{\newcommand{\Sgoth}{S\;\!\!\!\!\!/}
\newcommand{\bee}{\begin{equation}}
\newcommand{\ee}{\end{equation}}
\newcommand{\lam}{\lambda}
\newcommand{\ka}{\kappa}
\newcommand{\al}{\alpha}
\newcommand{\th}{\theta}
\newcommand{\fr}{\frac{1}{2}}
\newcommand{\fs}{\f{1}{2}}
\newcommand{\g}{\Gamma}
\newcommand{\br}{\biggr}
\newcommand{\bl}{\biggl}
\newcommand{\ra}{\rightarrow}
\newcommand{\mbint}{\frac{1}{2\pi i}\int_{c-\infty i}^{c+\infty i}}
\newcommand{\mbcint}{\frac{1}{2\pi i}\int_C}
\newcommand{\mboint}{\frac{1}{2\pi i}\int_{-\infty i}^{\infty i}}
\newcommand{\gtwid}{\raisebox{-.8ex}{\mbox{$\stackrel{\textstyle >}{\sim}$}}}
\newcommand{\ltwid}{\raisebox{-.8ex}{\mbox{$\stackrel{\textstyle <}{\sim}$}}}
\renewcommand{\topfraction}{0.9}
\renewcommand{\bottomfraction}{0.9}
\renewcommand{\textfraction}{0.05}
\newcommand{\mcol}{\multicolumn}
\date{}
\maketitle
\pagestyle{myheadings}
\markboth{\hfill \sc R. B.\ Paris  \hfill}
{\hfill \sc  Asymptotics of a tunnelling integral\hfill}
\begin{abstract}
We obtain an asymptotic evaluation of the integral
\[\int_{\sqrt{2n+1}}^\infty e^{-x^2} H_n^2(x)\,dx\]
for $n\ra\infty$, 
where $H_n(x)$ is the Hermite polynomial. This integral is used to determine the probability
for the quantum harmonic oscillator in the $n$th energy eigenstate to tunnel into the classically forbidden region.
Numerical results are given to illustrate the accuracy of the expansion.

\vspace{0.4cm}

\vspace{0.3cm}

\end{abstract}

\vspace{0.3cm}

\noindent $\,$\hrulefill $\,$

\vspace{0.2cm}

\begin{center}
{\bf 1. \  Introduction}
\end{center}
\setcounter{section}{1}
\setcounter{equation}{0}
\renewcommand{\theequation}{\arabic{section}.\arabic{equation}}
The quantum harmonic oscillator in dimensionless variables has the Hamiltonian given by ${\hat H}=\fs({\hat p}^2+{\hat x}^2)$, where ${\hat x}, {\hat p}$ are the position and momentum operators satisfying $[{\hat x}, {\hat p}]=i$. The normalised eigenstates are
\[\psi_n(x)=\pi^{-1/4}(2^n n!)^{-1/2} e^{-x^2/2} H_n(x),\]
where $H_n(x)$ is the $n$th Hermite polynomial. The corresponding probability densities $P_n(x)$ are then given by $|\psi_n(x)|^2$, that is
\bee\label{e11}
P_n(x)=a_n e^{-x^2} H_n^2(x),\qquad a_n=\frac{1}{\surd\pi 2^n n!}.
\ee

The classical turning points are situated at $x=\pm \sqrt{2n+1}$, so that the probability of tunnelling into the classically forbidden region is expressed as
\bee\label{e12}
P_{n,tun}=2a_nQ_n, \qquad Q_n=\int_{\sqrt{2n+1}}^\infty e^{-x^2} H_n^2(x)\,dx.
\ee
It has been pointed out recently by Jadczyk \cite{J} that the tunnelling probabilities $P_{n, tun}$ are rarely discussed in the literature on quantum mechanics. In \cite{J}, he has shown that $P_{n, tun}$ has the approximate representation for large values of $n$ given by
\bee\label{e13}
P_{n, tun}=\frac{1}{n^{1/3}}\left\{0.133975-\frac{0.0122518}{n^{2/3}}+O(n^{-1})\right\} \qquad (n\ra\infty).
\ee 

Our aim in this note is to obtain more precise information on the large-$n$ behaviour of $P_{n, tun}$. We carry this out by exploiting the uniform asymptotic expansion of the parabolic cylinder function $U(-n-\fs, 2^\fr x)$ given in \cite{T}, valid for $n\ra\infty$ and $x\geq \sqrt{2n+1}$ .
\vspace{0.6cm}

\begin{center}
{\bf 2. \ A representation of $Q_n$}
\end{center}
\setcounter{section}{2}
\setcounter{equation}{0}
\renewcommand{\theequation}{\arabic{section}.\arabic{equation}}
We set $\nu:=\sqrt{2n+1}$ and rescale the variable $x$ to write the integral $Q_n$ in (\ref{e12}) as
\[Q_n=\nu \int_1^\infty e^{-(\nu x)^2} H_n^2(\nu x)\,dx.\]
The integrand can be expressed in terms of the parabolic cylinder function $U(a,z)$ by \cite[Eq. (18.15.28)]{DLMF}
\[e^{-(\nu x)^2/2} H_n(\nu x)=2^{\frac{1}{4}\nu^2-\frac{1}{4}} U(-\fs\nu^2, 2^\fr\nu x).\]

The asymptotic expansion of this last function in terms of the Airy function Ai($z)$ and its derivative is given by \cite[Eq. (3.23)]{T}
\bee\label{e21}
U(-\fs\nu^2,2^\fr \nu x) =\frac{2^{\frac{1}{4}\nu^2+\frac{1}{4}}e^{\nu^2/4}\g(\fs+\fs\nu^2)}{\nu^{\fr\nu^2+\frac{1}{6}}}
\left(\frac{\zeta}{x^2-1}\right)^{\!1/4}\,\Upsilon_\nu(\zeta),
\ee
\bee\label{e21a}
\Upsilon_\nu(\zeta)=\mbox{Ai}\,(\nu^{4/3} \zeta)\,F_\nu(\zeta)+\nu^{-8/3} \mbox{Ai}'\,(\nu^{4/3}\zeta) \,G_\nu(\zeta)
\ee
as $\nu\ra+\infty$ uniformly valid for $x\geq 1$,
where
\bee\label{e21b}
\zeta^{3/2}=\f{3}{4}x \sqrt{x^2-1}-\f{3}{4}\,\mbox{arccosh} x\qquad (x\geq 1).
\ee
The functions $F_\nu(\zeta)$ and $G_\nu(\zeta)$ have the asymptotic expansions
\[F_\nu(\zeta)\sim\sum_{k=0}^\infty f_k(\zeta) \nu^{-2k}, \qquad G_\nu(\zeta)\sim\sum_{k=0}^\infty g_k(\zeta) \nu^{-2k},\]
where, from \cite[\S 3.2]{T}, we have
\[f_0(\zeta)=1,\quad f_1(\zeta)=\frac{1}{24}, \quad f_2(\zeta)=a_1(\zeta)+\frac{1}{576},\]
\[g_0(\zeta)=b_0(\zeta),\quad g_1(\zeta)=\frac{b_0(\zeta)}{24}, \]
with
\bee\label{e22b}
a_1(\zeta)=\frac{1}{1152}\bl\{\frac{145+249x^2-9x^4}{(x^2-1)^3}-\frac{7x(x^2-6)}{(x^2-1)^{3/2}\zeta^{3/2}}-\frac{455}{4\zeta^3}\br\},
\ee 
\bee\label{e22a}
b_0(\zeta)=-\frac{1}{2\zeta^\fr}\left\{\frac{x(x^2-6)}{12(x^2-1)^{3/2}}+\frac{5}{24\zeta^{3/2}}\right\}.
\ee

Then, from (\ref{e21}) and (\ref{e21a}) it follows that
\[Q_n=\frac{2^{2n+1} (n!)^2 e^{n+\fr}}{\nu^{2n+\frac{1}{3}}} \int_1^\infty\left(\frac{\zeta}{x^2-1}\right)^{\!1/2}
\Upsilon_\nu^2(\zeta)dx.\]
The quantity $\zeta$ is a monotonically increasing function for $x\geq 1$ and $\zeta=0$ when $x=1$. Furthermore
\[\frac{d\zeta}{dx}=\left(\frac{\zeta}{x^2-1}\right)^{-1/2},\]
so that we can introduce the new integration variable $\zeta$ to find
\bee\label{e22}
Q_n=\frac{2^{2n+1} (n!)^2 e^{n+\fr}}{\nu^{2n+\frac{1}{3}}} \int_0^\infty\phi(\zeta)\,\Upsilon_\nu^2(\zeta)\,
d\zeta,\qquad \phi(\zeta):=\frac{\zeta}{x^2-1}.
\ee

From (\ref{e21a}), we have
\[\Upsilon_\nu^2(\zeta)=\mbox{Ai}^2(\nu^{4/3}\zeta)\,F_\nu^2(\zeta)+2\nu^{-8/3}\mbox{Ai}\,(\nu^{4/3}\zeta) \mbox{Ai}'(\nu^{4/3}\zeta)\,F_\nu(\zeta) G_\nu(\zeta)\]
\[\hspace{8cm}+\nu^{-16/3} \mbox{Ai}'{}^2(\nu^{4/3}\zeta)\,G_\nu^2(\zeta)\]
\[=\mbox{Ai}^2(\nu^{4/3}\zeta)\bl(1+\frac{1}{12\nu^2}+\frac{1\!+\!1152f_2(\zeta)}{576\nu^4}+O(\nu^{-6})\br)\]
\[+\nu^{-8/3}b_0(\zeta)[\mbox{Ai}^2\,(\nu^{4/3}\zeta)]'\bl(1+\frac{1}{12\nu^2}+O(\nu^{-4})\br)
+\nu^{-16/3}b_0^2(\zeta) \mbox{Ai}'{}^2(\nu^{4/3}\zeta) (1+O(\nu^{-2})).
\]
Retaining terms up to and including $\nu^{-14/3}$, we therefore find
\[\int_0^\infty \phi(\zeta)\,\Upsilon_\nu^2(\zeta)\,d\zeta=
\int_0^\infty \phi(\zeta)\,\bl\{\mbox{Ai}^2(\nu^{4/3}\zeta)\bl(1+\frac{1}{12\nu^2}+\frac{1\!+\!1152f_2(\zeta)}{576\nu^4}\br)\hspace{2cm}\]
\bee\label{e23}
\hspace{4cm}+\nu^{-8/3} b_0(\zeta)\,[\mbox{Ai}^2\,(\nu^{4/3}\zeta)]'+O(\nu^{-16/3})\br\}d\zeta
\ee
as $\nu\ra\infty$.
\vspace{0.6cm}

\begin{center}
{\bf 3. \ The asymptotic expansion of two integrals}
\end{center}
\setcounter{section}{3}
\setcounter{equation}{0}
\renewcommand{\theequation}{\arabic{section}.\arabic{equation}}
To evaluate the integrals appearing in (\ref{e23}) we employ the modification of Watson's lemma given in Olver's book  \cite[p.~337]{O}. This is possible since Ai($x)$ possesses the asymptotic behaviour
\[\mbox{Ai}\,(x)\sim \fs \pi^{-1/2} x^{-1/4} \exp [-\f{2}{3}x^{3/2}]\qquad (x\ra +\infty).\]
Inversion of (\ref{e21b}) with the help of {\it Mathematica} yields\footnote{The function $\zeta(x)$ is real for $x>-1$, being given by
$(-\zeta)^{3/2}=\frac{3}{4}(\arccos x-x\sqrt{1-x^2})$ in the interval $-1<x\leq 1$. The finite radius of convergence results from the singularity in the inversion at $x=-1$, where $\zeta=-(\frac{3}{4}\pi)^{2/3}$.}
\[x=1+2^{-1/3}\zeta-\frac{2^{-2/3}\zeta^2}{10}+\frac{11\zeta^3}{700}-\frac{823\cdot 2^{-1/3}\zeta^4}{12600}+ \ldots \qquad (|\zeta|<\zeta_0),\]
where $\zeta_0=(\f{3}{4}\pi)^{2/3}$; see \cite[Eq.~(12.10.41)]{DLMF}. Then we have
\bee\label{e23a}
\phi(\zeta)=\sum_{m=0}^\infty\alpha_m \zeta^m \qquad (|\zeta|<\zeta_0), 
\ee
where
\[\alpha_0=2^{-2/3},\quad\alpha_1=-\frac{1}{5},\quad \alpha_2=\frac{2^{5/3}}{35},\quad \alpha_3=-\frac{2^{10/3}}{225},\quad\alpha_4=\frac{1548}{67375}, \ldots\, .\]

Now for a fixed positive integer $M$, and with $t=\nu^{4/3}\zeta$,
\[\int_0^\infty \phi(\zeta)\,\mbox{Ai}^2(\nu^{4/3}\zeta)\,d\zeta=\sum_{m=0}^{M-1}\frac{ \alpha_m}{\nu^{\frac{4}{3}m+\frac{4}{3}}} \int_0^\infty t^m \mbox{Ai}^2(t)\,dt+R_M,\]
where
\[R_M=\int_0^\infty \phi_M(\zeta)\,\mbox{Ai}^2(\nu^{4/3}\zeta)\,d\zeta,\qquad \phi_M(\zeta):=\phi(\zeta)-\sum_{m=0}^{M-1} \alpha_m \zeta^m.\] 
It follows that $\phi_M(\zeta)=O(\zeta^M)$ as $\zeta\ra 0+$ and it is known \cite{J} that $\phi(\zeta)\leq 2^{-2/3}$ for $\zeta\geq 0$. Hence we can find a constant $C$ such that the inequality $|\phi_M(\zeta)|<C\zeta^M$ holds when $\zeta\geq 0$. Then
\[|R_M|\leq C\int_0^\infty \zeta^M \mbox{Ai}^2(\nu^{4/3}\zeta)\,d\zeta=O(\nu^{-\frac{4}{3}M-\frac{4}{3}})\]
as $\nu\ra\infty$,
upon use of the evaluation \cite[Eq.~(12.11.15)]{DLMF}
\bee\label{e24b}
\int_0^\infty t^m \mbox{Ai}^2(t)\,dt=\frac{2m!\,12^{-\frac{1}{3}m-\frac{7}{6}}}{\surd\pi \g(\f{1}{3}m+\f{7}{6})}\qquad (m=0, 1, 2, \ldots ).
\ee
Thus we obtain the asymptotic expansion
\bee\label{e24}
\int_0^\infty \phi(\zeta) \mbox{Ai}^2(\nu^{4/3}\zeta)\,d\zeta\sim \frac{2}{\surd\pi} \sum_{m=0}^\infty \frac{\alpha_m m!\,\nu^{-4m/3-4/3}}{12^{\frac{1}{3}m+\frac{7}{6}} \g(\f{1}{3}m+\f{7}{6})} \qquad (\nu\ra\infty).
\ee

A similar reasoning can be applied to the second integral on the right-hand side in (\ref{e23}), when we observe that the function $-b_0(\zeta)$ in (\ref{e22a}) decreases monotonically from the value $9\cdot 2^{-2/3}/140$ to zero  for $\zeta\geq 0$ (we omit these details). With the help of {\it Mathematica}, the product $\phi(\zeta)b_0(\zeta)$ has the series expansion
\[\phi(\zeta)b_0(\zeta)=-\sum_{m=0}^\infty \beta_m \zeta^m\qquad (|\zeta|<\zeta_0),\]
where the first few coefficients are
\[\beta_0=\frac{9}{280} 2^{-1/3},\quad \beta_1=-\frac{179}{6300}2^{-2/3}, \quad \beta_2= \frac{28687}{2425500}, \quad
\beta_3=-\frac{750979}{78828750} 2^{-1/3}. \]
Then we find
upon using an integration by parts, combined with (\ref{e24b}) and the result Ai$(0)=3^{-2/3}/\g(\f{2}{3})$,
that
\bee\label{e25}
\int_0^\infty \phi(\zeta)\,b_0(\zeta)\,[\mbox{Ai}^2\,(\nu^{4/3}\zeta)]'\,d\zeta\sim
\frac{\beta_0 (3\nu)^{-4/3}}{\g^2(\f{2}{3})}+\frac{2}{\surd\pi}\sum_{m=1}^\infty \frac{\beta_m m!\, \nu^{-4m/3-4/3}}{12^{\frac{1}{3}m+\frac{5}{6}} \g(\f{1}{3}m+\f{5}{6})}
\ee
as $\nu\ra\infty$.
\vspace{0.6cm}

\begin{center}
{\bf 4. \ The asymptotic calculation of $P_{n, tun}$ for $n\ra\infty$}
\end{center}
\setcounter{section}{4}
\setcounter{equation}{0}
\renewcommand{\theequation}{\arabic{section}.\arabic{equation}}
The function $-a_1(\zeta)$ appearing in (\ref{e22b}) is found to be monotonically decreasing for $\zeta\geq 0$ with the series expansion
\[-a_1(\zeta)=\frac{249}{28800}-\frac{6849}{616000} 2^{-1/3}\zeta+\frac{737}{65000} 2^{-2/3}\zeta^2 - \cdots\qquad(|\zeta|<\zeta_0).\]
Then the O($\nu^{-4})$ term multiplying $\mbox{Ai}^2(\nu^{4/3}\zeta)$ in (\ref{e23}) has the expansion
\[-\frac{1}{576}\phi(\zeta)(1+1152f_2(\zeta))=\frac{29}{2400} 2^{-2/3}-\frac{25013}{1848000} \zeta+ \cdots \qquad (|\zeta|<\zeta_0).\]

From (\ref{e23}), (\ref{e24}) and (\ref{e25}) we therefore find
\[\nu^{4/3}\int_0^\infty \phi(\zeta) \Upsilon_\nu^2(\zeta)\,d\zeta\hspace{8cm}\]
\[=\frac{2}{\surd\pi}\bl(1+\frac{1}{12\nu^2}-\frac{29}{2400\nu^4}\br)
\bl(\frac{1}{12^{7/6} \g(\f{7}{6})}+\frac{2\alpha_1 \nu^{-4/3}}{12^{3/2} \surd\pi}+\frac{2\alpha_2 \nu^{-8/3}}{12^{11/6} \g(\f{11}{6})}
+\frac{6\alpha_3\nu^{-4}}{12^{13/6} \g(\f{13}{6})}\br)\]
\[+ \nu^{-8/3}\bl(1+\frac{1}{12\nu^2}\br) \bl(\frac{3^{-4/3}\beta_0}{\g^2(\f{2}{3})}+\frac{2\beta_1 \nu^{-4/3}}{\surd\pi\,12^{7/6} \g(\f{7}{6})}\br)+O(\nu^{-16/3}).\]
Hence, from (\ref{e12}) and (\ref{e22}), we obtain
\[P_{n, tun}=\frac{2^{n+2} (n!) e^{n+\fr}}{\surd\pi\,\nu^{2n+5/3}}\hspace{8cm}\]
\[\times\left\{\bl(1+\frac{1}{12\nu^2}-\frac{29}{2400\nu^4}\br) \bl\{\frac{6^{-2/3}}{\g^2(\f{1}{3})} -\frac{\nu^{-4/3}}{30\pi \surd 3}+\frac{4}{525}\,\frac{6^{-1/3} \nu^{-8/3}}{ \g^2(\f{2}{3})}-\frac{16}{525}\,\frac{6^{-2/3} \nu^{-4}}{\g^2(\f{1}{3})}\br\}\right.\]
\bee\label{e41}
\left.+\nu^{-8/3}\bl(1+\frac{1}{12\nu^2}\br)\bl\{\frac{3}{280}\,\frac{6^{-1/3}}{\g^2(\f{2}{3})}-\frac{179}{6300}\,\frac{6^{-2/3} \nu^{-4/3}}{\g^2(\f{1}{3})}\br\}+O(\nu^{-16/3})\right\}
\ee
upon use of the properties of the gamma function to simplify the coefficients.

On noting that
\[\frac{2^{n+2} n! e^{n+\fr}}{\surd\pi\,\nu^{2n+5/3}}=\frac{2^{5/3}}{n^{1/3}}\bl(1-\frac{5}{12\nu^2}-\frac{23}{288\nu^4}+O(\nu^{-6})\br)\qquad (n\ra\infty),\]
we then finally obtain from (\ref{e41}) that
\bee\label{e42}
P_{n, tun}=\frac{2^{5/3}}{n^{1/3}}\bl\{\bl(1-\frac{1}{3\nu^2}\br)\bl(\frac{6^{-2/3}}{\g^2(\f{1}{3})} -\frac{\nu^{-4/3}}{30\pi \surd 3}+\frac{11}{600}\,\frac{6^{-1/3} \nu^{-8/3}}{ \g^2(\f{2}{3})}\br)-\frac{167}{900}\,\frac{6^{-2/3}\nu^{-4}}{\g^2(\f{1}{3})}+O(\nu^{-16/3})\br\}
\ee
as $n\ra\infty$. Insertion of numerical values for the coefficients yields the alternative form
\[P_{n, tun}=\frac{1}{n^{1/3}}\left\{\bl(1-\frac{1}{3\nu^{2}}\br)\bl(0.1339750-\frac{0.0194484}{\nu^{4/3}}-\frac{0.0174687}{\nu^{8/3}}\br)\hspace{2cm}\right.\]
\[\hspace{7cm}\left.-\frac{0.0248598}{\nu^{4}}+O(\nu^{-16/3})\right\},\quad \nu=\sqrt{2n+1}.\]
The first two leading terms in this expansion can be seen to agree with the approximation stated in (\ref{e13}).

We conclude by presenting the results of numerical calculations to illustrate the accuracy of the asymptotic expression in (\ref{e42}).
In Table 1 we show  the values of the absolute relative error in the computation of $P_{n, tun}$ using (\ref{e42}) for different $n$. 
\begin{table}[bh]
\caption{\footnotesize{Values of the absolute relative error in the computation of $P_{n, tun}$ from (\ref{e42}).}}
\begin{center}
\begin{tabular}{|l|l|l||l|l|l|}
\hline
&&&&&\\[-0.25cm]
\mcol{1}{|c|}{$n$} & \mcol{1}{c|}{$P_{n,tun}$} & \mcol{1}{c||}{Error} & \mcol{1}{c|}{$n$}  & \mcol{1}{c|}{$P_{n, tun}$} & \mcol{1}{c|}{Error}\\
[.1cm]\hline
&&&&&\\[-0.3cm]
10 & 0.0601438 & $1.323\times 10^{-5}$ & 200 & 0.0228302 & $6.975\times 10^{-9}$\\
20 & 0.0483977 & $2.528\times 10^{-6}$ & 400 & 0.0181454 & $1.130\times 10^{-9}$\\
50 & 0.0360132 & $2.534\times 10^{-7}$ & 500 & 0.0168499 & $6.276\times 10^{-10}$\\
100& 0.0286973 & $4.250\times 10^{-8}$ & 800 & 0.0144138 & $1.815\times 10^{-10}$\\
[.2cm]\hline
\end{tabular}
\end{center}
\end{table}

\end{document}